\documentclass[a4paper,10pt]{article}
\usepackage{authblk}

\usepackage[ngerman,english]{babel}
\usepackage[latin1]{inputenc}
\usepackage{csquotes}
\usepackage{amsfonts,amsmath,amsthm}
\usepackage{empheq}
\usepackage[titletoc,title]{appendix}
\usepackage[backend=bibtex8,doi=false,eprint=false,firstinits=true,isbn=false,style=numeric-comp,maxnames=99]{biblatex}
\makeatletter
\def\blx@maxline{77}
\makeatother

\bibliography{bibl.bib}
\AtEveryBibitem{\clearlist{language}}

\usepackage{cases}
\usepackage{mathabx}
\usepackage{bbm}
\usepackage{xfrac}
\usepackage{fancyhdr}
\usepackage{color}
\usepackage[colorinlistoftodos,textsize=small,backgroundcolor=white,bordercolor=blue,linecolor=blue,disable]{todonotes}
\usepackage[colorlinks]{hyperref}
\definecolor{blue75}{rgb}{0,0,.75}
\definecolor{green75}{rgb}{0,.75,0}
\hypersetup{colorlinks=true, urlcolor=blue75,linkcolor=blue75,citecolor=green75,pdfstartview=FitB,bookmarksopen=true,bookmarksopenlevel=1}
\usepackage[a4paper, left=2.5cm, right=2.5cm, top=2.5cm,bottom=2cm]{geometry}
\usepackage{constants}
\newcommand{\parenthezises}[1]{\arabic{#1}}
\newconstantfamily{C}{
symbol=C,
format=\parenthezises,
reset={section}
}
\newconstantfamily{M}{
symbol=M,
format=\parenthezises,
reset={section}
}
\newconstantfamily{B}{
symbol=B,
format=\parenthezises,
reset={section}
}
\newconstantfamily{xi}{
symbol=\xi,
format=\parenthezises,
reset={section}
}
\usepackage{enumerate}
\usepackage{graphicx}
\graphicspath{{images/} }
\usepackage{wrapfig}
\usepackage{figbib}
\allowdisplaybreaks
\usepackage[capitalise]{cleveref}

\crefdefaultlabelformat{{\it #2#1#3}}

\crefname{equation}{}{}

\crefname{enumi}{}{}
\creflabelformat{enumi}{{(#2#1#3)}}

\crefname{section}{{\it Section}}{{\it Sections}}
\crefname{subsection}{{\it Subsection}}{{\it Subsections}}
\crefname{subsubsection}{{\it Paragraph}}{{\it Paragraphs}}

\newtheorem{Theorem}{Theorem}[section]
\crefname{Theorem}{{\it Theorem}}{{\it Theorems}}
\newtheorem{Definition}[Theorem]{Definition}
\crefname{Definition}{{\it Definition}}{{\it Definitions}}

\crefname{Lemma}{{\it Lemma}}{{\it Lemmas}}
\newtheorem{Proposition}[Theorem]{Proposition}
\crefname{Proposition}{{\it Proposition}}{{\it Propositions}}
\newtheorem{Assumptions}[Theorem]{Assumptions}
\crefname{Assumptions}{{\it Assumptions}}{{\it Assumptions}}

\theoremstyle{definition}
\newtheorem{Remark}[Theorem]{Remark}
\crefname{Remark}{{\it Remark}}{{\it Remarks}}

\crefname{Notation}{{\it Notation}}{{\it Notations}}

\crefname{Example}{{\it Example}}{{\it Examples}}

\renewbibmacro*{doi+eprint+url}{%
    \printfield{doi}%
    \newunit\newblock%
    \iftoggle{bbx:eprint}{%
        \usebibmacro{eprint}%
    }{}%
    \newunit\newblock%
    \iffieldundef{doi}{%
        \usebibmacro{eprint}}%
        {}%
    }
\begin{document}

\newcommand{\cb}[1]{{\color{blue}{#1}}}
 \newcommand{\red}[1]{\textcolor{red}{#1}}
\newcommand{\cmg}[1]{{\color{magenta}{#1}}}
\newcommand{\cgr}[1]{{#1}}
\newcommand{\D}{\mathbb{D}}
\newcommand{\E}{\mathbb{E}}
\newcommand{\T}{\mathbb{T}}
\newcommand{\PP}{\mathbb{P}}
\newcommand{\R}{\mathbb{R}}
\newcommand{\N}{\mathbb{N}}
\newcommand{\F}{\mathbb{F}}
\newcommand{\V}{\mathbb{V}}
\newcommand{\ve}{\varepsilon}

\def\diam{\operatorname{diam}}
\def\dist{\operatorname{dist}}
\def\diver{\operatorname{div}}
\def\ess{\operatorname{ess}}
\def\inner{\operatorname{int}}
\def\osc{\operatorname{osc}}
\def\sign{\operatorname{sign}}
\def\supp{\operatorname{supp}}
\newcommand{\BMO}{BMO(\Omega)}
\newcommand{\LOne}{L^{1}(\Omega)}
\newcommand{\LOnen}{(L^{1}(\Omega))^d}
\newcommand{\LTwo}{L^{2}(\Omega)}
\newcommand{\Lq}{L^{q}(\Omega)}
\newcommand{\Lp}{L^{2}(\Omega)}
\newcommand{\Lpn}{(L^{2}(\Omega))^d}
\newcommand{\LInf}{L^{\infty}(\Omega)}
\newcommand{\HOneO}{H^{1,0}(\Omega)}
\newcommand{\HTwoO}{H^{2,0}(\Omega)}
\newcommand{\HOne}{H^{1}(\Omega)}
\newcommand{\HTwo}{H^{2}(\Omega)}
\newcommand{\HmOne}{H^{-1}(\Omega)}
\newcommand{\HmTwo}{H^{-2}(\Omega)}

\newcommand{\LlogL}{L\log L(\Omega)}

\def\avint{\mathop{\,\rlap{-}\!\!\int}\nolimits} 

\newcommand{\om}{\omega }
\newcommand{\Om}{\Omega }

\newtheorem{proofpart}{Step}
\makeatletter
\@addtoreset{proofpart}{Theorem}
\makeatother
\numberwithin{equation}{section}
\title{A novel derivation of rigorous macroscopic limits from a micro-meso description of signal-triggered cell migration in fibrous environments}
	
\author{Anna Zhigun\thanks{School of Mathematics and Physics, Queen's University Belfast, University Road, Belfast BT7 1NN, Northern Ireland, UK, \href{mailto:A.Zhigun@qub.ac.uk}{A.Zhigun@qub.ac.uk}} \ and
    Christina Surulescu\thanks{Felix-Klein-Zentrum für Mathematik, Technische Universität Kaiserslautern, Paul-Ehrlich-Str. 31, 67663 Kaiserslautern, Germany, \href{mailto:surulescu@mathematik.uni-kl.de}{surulescu@mathematik.uni-kl.de}} 
   }

\date{}
\maketitle
\begin{abstract}
In this work we upscale a prototypical kinetic transport equation which models a cell population moving in a fibrous environment with a chemo- or haptotactic signal influencing both the direction and the magnitude of the cell velocity.  The presented  approach to scaling does  not rely on orthogonality and  treats  parabolic and hyperbolic scalings in a unified manner. It is shown that the steps of the  formal limit procedures are mirrored by rigorous operations with finite measures provided that the measure-valued position-direction fiber  distribution enjoys some spacial continuity. 
\end{abstract}

\section{Introduction}\label{sec:intro}
\subsection{{Biological motivation and modelling aspects}}
Cell migration is a highly complex biological process involving a multitude of mechanisms and being influenced by manifold chemical and physical components of the extracellular environment. Cell motility is decisive for a plethora of physiological processes {such as} wound healing, regeneration of tissues, embryonic development, tumour growth and metastasis. Most of these involve interactions of the respective cells (e.g. fibroblasts, endothelial cells, chondrocytes, osteoblasts, tumour cells) {with fibrous 
parts of the surrounding} tissue, which they are able to modify in various ways, but also use as support and guidance for migration, proliferation, and even survival. Most cells adapt dynamically to changes in soluble (chemoattractants and/or -repellents) and insoluble (tissue) components of their environment, thus exhibiting a tactic behavior. Thereby they can perceive space-time variations of the respective signals, thus follow concentration/density gradients. Such motile behavior mediated by tissue and chemical signals is commonly termed haptotaxis and chemotaxis, respectively, and has {been} the object of numerous mathematical works concerned with modeling these phenomena and performing analysis and/or numerics for the obtained systems of differential equations. We refer, e.g. to \cite{BeBeTaWi,HiPa08,Horstmann,KSSSL,Painter2019} for reviews of taxis models from various perspectives. 

One way to deduce such reaction-diffusion-taxis equations (RDTEs) is to consider a multiscale approach  {which starts from the microscopic scale of single cell behaviour (characterising cell trajectories and possibly also the so-called activity variables  \cite{bellom3})} and obtains one or several (depending on the number of cell (sub)populations involved) corresponding kinetic transport equations (KTEs){. Each of these equations is stated for a cell (sub)population distribution density} depending on time, position, velocity, and the mentioned activity variables. A {subsequent} appropriate upscaling leads from this mesoscopic description to effective RDTEs on the macroscopic scale of cell population(s). {Equations obtained in this manner inherit important low-level information from the original KTEs allowing for a more detailed and accurate modelling.} Models for cell dispersal which {were obtained in} this framework have been proposed, e.g. in  \cite{Chalub2004,Erban2005,HILLEN2002,Othmer2000,Perthame2018,Perthame2020,Xue2009} for various scenarios of bacteria performing chemotaxis, some of these works also  providing rigorous macroscopic limits. Still in the context of chemotactic behaviour, flux-limited macroscopic RDTEs have been obtained from KTEs e.g. in \cite{BBNS2010} by appropriate macroscopic limits; those systems have, among other advantages, the essential ability to enable finite propagation speed of the cell population performing the respective kind of diffusion and/or chemotaxis.

The dynamics of cell migration in fibrous tissues is more complex, as the cells perform several types of taxis, of which haptotaxis is essential \cite{carter}. Moreover, some cancer cells like glioma { (generating the most common type of primary brain tumors)} are able
to exploit the underlying tissue anisotropy in order to enhance their motility \cite{Giese1996}, thus the orientation of tissue fibers is relevant and should be included in the modelling.  First models employing KTEs and
addressing this issue in the framework of cell migration were proposed  in \cite{Chauviere2007,Hillen2006,HillenPainter,Painter2013}. Of these, in the former study a macroscopic {RDTE} involving
chemo- and haptotaxis was obtained, whereas in the others  some appropriate scaling led to  {RDTE}s with myopic diffusion explicitly involving the space-dependent distribution of tissue fibers. Works \cite{KELKEL2012,Lorenz2014} extended the common KTE modelling framework  to include both types
of taxis (via adequate kernels) in the terms characterizing velocity reorientations.
More recent multiscale models characterizing chemo- and/or haptotactic cancer cell migration in tissue networks and featuring myopic dffusion were proposed in \cite{CS20,Corbin2018,CEKNSSW,EHS,Engwer2,Engwer,Hunt2016,Kumar20}, all of which performed {formal} upscalings from the micro-meso{-} to the {macroscale}. 

In this work we deduce, first formally and then rigorously,  macroscopic diffusion-taxis equations (DTEs) for cell migration in an anisotropic tissue.  
{Following the multiscale modelling approach, we first construct a (mesoscopic) conservative linear KTE   involving
 transport with respect  to velocity and then}
 perform parabolic (provided that the mean fiber orientation vanishes) and hyperbolic upscalings.
The equations obtained in the macroscopic limit are accordingly  diffusion- or  drift-dominated. Recently,  \cite{CEKNSSW}  shortly reviewed the hitherto available alternative ways to include cell level environmental influences in a {KTE-based} modelling framework leading  to taxis terms in the macroscopic limit. Of those, the present note addresses the one which accounts for biochemical and/or biophysical effects translated into cell stress and forces acting on the cells. 
In the corresponding KTE this results in a non-zero transport term with respect to velocity.
Such  approach was also considered, e.g. in \cite{Chauviere2007,CEKNSSW,DKSS20}, the latter two actually combining it with the involvement of some further cell activity variables. Our {model}  is closer to \cite{DKSS20},  {yet accounts for less biological complexity. In particular, it does not include any activity  variables. Nevertheless,  cells considered here  can {\it actively} change their velocity due to the presence of a stimulus. The aim of our   work is to perform formal and rigorous upscalings  for the model, presenting  them in a transparent manner. }
\subsection{Mathematical  aspects}
\subsubsection{Macroscopic approximations}
{Modelling with mesoscopic KTEs and subsequently upscaling them to get macroscopic RDTEs} has proved to be a very effective scheme in the context of  population migration in heterogeneous surroundings{. However,} 
the calculations {are} usually performed only formally and lacked rigorous justification. 
Overall, those settings which are very accurate from the modelling point of view are often  particularly  difficult to handle rigorously. One reason here is that since a rigorous limit procedure amounts to a proof of existence of solutions to the resulting macroscopic  model, it is not surprising that it can be very difficult,  if not impossible, to carry out if that model is nontrivial. An example of this is  the study \cite{WinSur2017} of a one-dimensional haptotaxis system with a  degenerate myopic diffusion arising from a potential tissue heterogeneity which was formally obtained through upscaling in \cite{EngHuntSur}. In this model, the main challenges are a strong coupling of the equations and the strongly degenerate diffusion. 
 Yet another illustrative example is a system describing   particle motion in a fluid for which rigorous parabolic and hyperbolic scalings were conducted in  \cite{GoudJabVass1} and \cite{GoudJabVass2}, respectively. There the challenges on the macroscopic level are due to a strong coupling and the fact that one of the equations is an incompressible Navier-Stokes equation. 
 
In this note we avoid the coupling issue and consider a single mesoscopic linear KTE for a single cell population moving in a heterogeneous external environment under the influence of an external signal. The biophysical/biochemical force in our KTE has a similar form  to the drag Stokes force used in \cite{GoudJabVass1,GoudJabVass2}. It is proportional to the difference between the cell velocity and  an external velocity field. In our case the field stems from a stimulus, not a fluid. Furthermore, unlike that model, not a diffusion with respect to velocity, but, rather, a turning operator is incorporated into the KTE. The latter describes the chaotic interactions with a heterogeneous extracellular environment, leading to a possibly degenerate myopic diffusion as in \cite{HillenPainter,WinSur2017}.
To summarise, our equation can be categorised as a Vlasov-Fokker-Planck equation into which a linear turning operator has been incorporated. From the mathematical perspective, it can also be viewed as a variant of the linear 
Boltzmann-Maxwell equation, a linear kinetic equation with an external field. Such an equation describes a gas of charged particles moving under the influence of an external field  through an unchanging background of another type of particles. However, it should be stressed that the interactions between cells and, e.g. tissue fibers are very different from collisions between two kinds of physical particles. In particular, the directional  distribution of the fibers is by no means a Maxwellian and can  actually be very unsmooth. In this work we allow it to be a finite measure, such  as, e.g. a Dirac delta function.

Establishing convergence of a sequence of linear KTEs such as ours
to, e.g. a parabolic RDE can be complicated, as these equations are of a very   different nature. Indeed, a KTE is an integro-differential equation that includes both a divergence operator, typically with respect to the time and space variables (in the  transport part),  and integration (in the turning operator) with respect to the velocity variable which needs to be scaled out.  
In our case due to the external forces the divergence with respect to the velocity variable  is also present in the transport. This  makes the upscaling even more challenging.

Currently, the most commonly applied {up}scaling procedures in the context of modelling 
{cell} movement are: the parabolic scaling, the hyperbolic scaling, and the moment closure method. Since our approach relies on {rescaling}, we do not consider the latter method and concentrate on the first two instead. {As usual for such scalings, a small scaling parameter $\varepsilon$ is introduced, the KTE rescaled, and the aim is then to approximate the solutions as $\varepsilon$ is sent to zero. The approximations are typically sought in the form of Hilbert or Chapman-Enskog expansions. 
The traditional approach to both parabolic and hyperbolic limits relies on some Hilbert space structure already for the formal scalings, see e.g. the derivations presented in \cite{HillenPainter,Othmer2000}.
 In particular, the Chapman-Enskog expansion that is typically used in the hyperbolic case   is based on the assumption that the actual limit function, i.e. the zero order term in the  expansion, is orthogonal to the first order correction.}  However, this assumption cannot be justified in such models as the one we consider here. {Our} calculations show {that} the resulting leading zero order term and the first order correction are not mutually orthogonal in general. The  first order correction has a rather complicated form, carrying the dependence  on the gradient of {an external} macroscopic quantity.
 
 A technique that has been extensively used when  dealing with rigorous scalings of linear KTEs is the semigroup theory  \cite{MikaBanasiak}.  However, the presence of the derivative with respect to the velocity variable turns out to be problematic already in very basic cases, such as the hyperbolically scaled one-dimensional linear KTE (10.1.9) from \cite[Chapter 10]{MikaBanasiak}. Unlike our model, there the velocity field is taken to be constant, so that it does not depend on the velocity variable, and a regular Maxwellian distribution is used in the collision term.  Already in this simplified situation it turns out necessary to work in the phase space which is a weighted $L^1$ space in order to pass rigorously to a macroscopic limit. In our case, we are forced to work in even rougher spaces of Radon measures.

In this note we present a novel and unified approach towards both parabolic and hyperbolic scalings, which  avoids the necessity to work in Hilbert spaces and does not  use semigroups. Its key component is a differential equation connecting moments of order zero and two of the KTE solutions, which at  first glance  resembles a damped wave equation with a transport term. Differentiating this equation a sufficient number of times with respect to the scaling parameter and passing to the limit allows to obtain the coefficients of the Hilbert expansion of any order. 
It turns out that thanks to the linearity of the KTE and the mass preservation, under rather general assumptions  the derivation steps can then be mirrored by operations with finite measures leading to a rigorous scaling. Our work extends a result obtained  in \cite{HHW10} where a parabolic scaling was performed  in the special case of a (non-evolving) space-homogeneous fiber distribution and for a KTE without velocity derivative.

\subsubsection{Mesoscopic approximation}
Often, when considering e.g., Hilbert expansions of the form
\cgr{$c^{\ve}=c^{0}+\ve c^{0}_1+O(\varepsilon^2)$}
for the solution of the original KTE, one would like to combine {the
leading} order term $c^0$ and  first order  correction $\ve c^{0}_1$ in order to get a mesoscopic approximation. 
The basic  straightforward approach  which consists of directly adding them together is known to have a serious drawback: the resulting function is not necessarily nonnegative. One  way to avoid this issue is to use a nonlinear Hilbert expansion as was done, e.g. in \cite{CGCT2005} for a simpler equation and in space dimension one. In this work we develop an alternative approach which requires  dealing with a truncated version of the original mesoscopic equation. This new PDE is a transport equation in terms of a mesoscopic first order approximation. It preserves both positivity and the total mass and can be solved numerically after the macroscopic approximations of zero and first orders are obtained. 
\\

The remainder of the paper is structured as follows. To begin with, we introduce in  \cref{SecKTE} a prototypical KTE describing cell movement on the mesolevel. For this model,  we consider in \cref{SecScaling} parabolic and hyperbolic scalings of time and space and  formally derive the corresponding limits and first-order corrections,  also developing  a  mesoscopic first order approximation which preserves both positivity and the total mass. The rigorous scaling is done in \cref{SecRig}. In the closing \cref{SecDis} we discuss our findings.

\section{A KTE modeling framework}\label{SecKTE}
{In this Section we describe the KTE modelling framework which is the starting point for our approach. More precisely, we} consider the migration of cancer cells in a fibrous, anisotropic tissue under the influence of some extracellular signal. 
The description of single cell behaviour involves the position and velocity dynamics which allow to reconstruct cell trajectories. Without loss of generality we assume {the cell speed magnitude to be less than one}
and consider the velocity space{
\begin{align}
 V=B_{1}(0):=\{v=s\theta:\ s\in[0,1),\ \theta\in S_1(0)\},\quad  S_1(0):=\{\theta\in\R^n: |\theta|=1\}.
\end{align}}
For any quantity $u$ depending on $v$ we will denote
\begin{align}
 \overline u:=\int_{B_{1}(0)}u(v)\,dv.
\end{align}
{In order to  capture { with some accuracy} the evolution of a large cell population, we start our modelling on the microscopic scale. Disregarding for the moment the chaotic interactions with tissue, we assume the movement of each cell to be  described by the following system of ordinary differential equations} { corresponding to Newton's second law:}
\begin{subequations}\label{MicroD}
\begin{align}
&\frac{dx}{dt}=v, \\
&\frac{dv}{dt}={S(t,x,v)},\label{eq:microdyn}
\end{align}
\end{subequations}
where
\begin{align}
 &{S(t,x,v)=-a(v-v_*(t,x))},\qquad  a>0,\ v_*{(t,x)}\in B_{1}(0){ \text{ for all }t>0,\ x\in\R^n}.\label{S}
\end{align}
{Here $t$, $x$, $v$, and $S$ are time, position in space, velocity, and acceleration of a cell, respectively. Our choice of $S$ is similar to the one in \cite{DKSS20}.} 
\cgr{It reflects the cell tendency to redirect its velocity in order to realign with a certain preferred direction $v_*$ and to decelerate if such a direction is missing. } 
Most of the previous constructions assumed zero acceleration or such that  vanishes at least in the radial direction, like, e.g. in \cite{CEKNSSW}, or does not depend on $v$, like in \cite{Chauviere2007}. For \cref{S} this is clearly not the case, so that in our model a cell can change not only the direction, but also 
 its speed and the change depends on the velocity itself. As we will see later in \cref{SecScaling}, this makes calculations more involved. 
%

\noindent
As in \cite{Chauviere2007,DKSS20}, we model $v_*$ in such a way as to account for
the effect of chemo- or haptotaxis. 
The model proposed in \cite{Chauviere2007} considered for the cell velocity dynamics a 'chemotactic force' depending on the gradient of some given chemical profile, hence also allowing for changes in the cell speed. The setting in \cite{DKSS20} combines such 'tactic forces' (repellent chemotaxis and haptotaxis) with a repellent force caused by large gradients of macroscopic cell density, all these influences contributing to modifications of direction as well as speed of the migrating cells.\\
 {Aiming to describe the  effect of the (macroscopic) tissue density or, alternatively, of some (possibly different) signal 
 on the cell  velocity reorientation, here} 
 we choose  $v_*$ in the form
\begin{align}
v_*=\F\frac{\nabla_x Q}{1+|\nabla_x Q|},\label{vstar}
\end{align}
where $Q=Q(t,x)$ is the concentration of a chemoattractant or the  macroscopic density of tissue {(one could also consider both,  introducing  yet another macroscopic quantity, $\mathbb A(t,x)$ for that purpose)} and $\F=\F(x)$ is a tensor such that
\begin{align*}
 \|\F(x)\|_{2}\leq 1\qquad\text{for all }x\in\R^n.
\end{align*}
 As in some previous works  \cite{CEKNSSW,DKSS20},
 the tensor $\F$ is supposed to model, e.g. biomechanical cell stress or simply have a deviatoric effect due to the heterogeneity of the environment. 
 {For instance, when modelling the glioma spread,  one could choose (similarly to \cite{DKSS20}), { for each brain voxel centered at $x$, the matrix} $\mathbb F{ (x)}$ to be the water diffusion tensor obtained by diffusion tensor imaging (DTI), a variant of MRI which is a standard noninvasive diagnostic tool. There is abundant evidence \cite{Giese1996} that migrating glioma cells follow  the highly anisotropic brain structure, particularly that of white matter.
 One could also choose the tensor { $\mathbb F$} to be  velocity-dependent, e.g. taking   $$\F:=\frac{1}{a}\left (|v|^2\mathbb I_n-v\otimes v\right),$$ as it was done in \cite{CEKNSSW}. Here and in what follows  $\mathbb I_n$ denotes the identity matrix and $\otimes$ stands for the tensor product. We refer to that work for further details.
 In any case, the eigenvalues and eigenvectors of $\F$ are used to encode relevant information about local anisotropy and diffusivity, thus helping to reconstruct the whole structure of brain tissue and opening the way for predictions about the extent of tumour spread. For more details on this approach we refer to \cite{CEKNSSW,Engwer2,EHS,Painter2013}.
 }

We model the fibrous { extracellular} environment by way of the orientational distribution of tissue fibers,  upon taking
\begin{align}
  &q=q\left(x,\hat v\right),\quad \hat v=\frac{v}{|v|},\label{qsign}\\
 &q\geq0,\quad  
 \int_{S_1(0)}q(\theta )\,d\theta =1,\label{qdistr}
\end{align}
hence 
\begin{align}
 \int_{B_{1}(0)}q({ v})\,dv=\frac{1}{n}.
\end{align}
For later purposes we introduce the moments 
\begin{align}
 &\E[q]:=\int_{S_1(0)}\theta q(\theta )\,d\theta,\\
 &\D[q]:=\int_{S_1(0)}\theta \otimes \theta q(\theta )\,d\theta,\\
 &\T[q]:=\int_{S_1(0)}\theta \otimes \theta \otimes \theta q(\theta )\,d\theta.\label{TenT}
\end{align}
Here $\E[q]$ can be interpreted as the average orientation of tissue fibers, while $\D[q]$ {is the auto-correlation matrix and thus related to the} variance-covariance matrix $$\mathbb V[q]:=\int_{S_1(0)}(\theta-\E[q]) \otimes (\theta -\E[q] )\,q(\theta)\ d\theta,$$ so that $$\D[q]=\mathbb V[q]+\E[q]\otimes \E[q]. $$
As usual, when performing a parabolic scaling, we assume that
$$ \E[q]\equiv 0,$$
which implies that 
\begin{align}
\mathbb V[q]\equiv\D[q].\label{VD}\end{align}

\noindent
\cgr{Finally, if $\E[q]=0$, then the third moment, $\T[q]$ accounts for the skewness  of distribution $q$, i.e. for its  asymmetry about its mean. }

\noindent
{Now we have all ingredients necessary to construct a model on the mesoscale. On this level we model the probability density function $c=c(t,x,v)$ for {cell} position and velocity  over time. 
 {Our equation for $c$ reads:} 
\begin{align}
 \nabla_{(t,x,v)}\cdot((1,v,{S(t,x,v)})c)=\partial_t c+\nabla_x\cdot (vc)-a\nabla_v\cdot((v-v_*)c)
 =&nq\overline{c}-c.\label{meso}
\end{align}
{The left-hand side of \cref{meso} describes the {mass-conserving} transport 
{along the trajectories generated by}
\cref{MicroD}. The right-hand side captures  chaotic interactions with tissue fibers. This is done by means of incorporating  a turning operator with a}
 turning kernel $nq$ and a constant turning rate (which for simplicity we set to 1). Thus, as in many previous works (see e.g. \cite{CS20,Corbin2018,CEKNSSW,EHS,Engwer2,Engwer,Hillen2006, Hunt2016,Kumar20,Painter2013}) we assume that the cells adapt their respective direction of motion to the local orientation of tissue fibers. Moreover, still in line with previous works, we consider the cell {population} to be compactly supported on the velocity space, thus introduce the boundary condition
\begin{align}
 c=0\qquad\text{for }v\in\partial B_{1}\label{bc}.
\end{align}


\section{Formal upscaling}\label{SecScaling}
In this section  we derive a formal macroscopic limit for the suitably rescaled KTE \cref{meso}. A rigorous argument requires adequate assumptions on the model parameters and is based on calculations with Radon measures. We postpone this issue to  \cref{SecRig}.

In order to set the frame we introduce some macroscopic time and space scales: for $\ve\in(0,1]$ let 
\begin{align}
 &\hat{t}=\ve^{\kappa}t,\qquad\kappa\in\{1,2\},\\ &\hat{x}=\ve x.
\end{align}
%
Rescaling \cref{meso} and \cref{vstar} and dropping the hats leads to
\begin{align}
 &\ve^{\kappa}\partial_t c^{\ve}+\ve \nabla_x\cdot (vc^{\ve})-a\nabla_v\cdot((v-v_*^{\ve})c^{\ve})=nq\overline{c^{\ve}}-c^{\ve},\label{transpce}
\end{align}
where 
\begin{align}
 v_*^{\ve}={\F}\frac{\ve \nabla_x Q}{1+\ve|\nabla_x Q|}=\ve {\F}\nabla_x Q-\ve^2 {\F}\nabla_x Q|\nabla_x Q|+O\left(\ve^3\right).
 \label{vstare}
\end{align}

\subsection{Equations involving moments}\label{SubSMom}
To begin with, we integrate \cref{transpce} by parts with respect to $v$ over $B_{1}(0)$ and divide by $\ve^{\kappa}$ in order to obtain an equation which connects the moments of orders zero and one:
\begin{align}
 \partial_t \overline{c^{\ve}}+\ve^{1-\kappa}\nabla_x\cdot \overline{vc^{\ve}}=0.\label{mom0}
\end{align}
Thereby we used the fact that $c^{\ve}$ vanishes on the boundary. 
Next, we multiply \cref{transpce} by $v$ and  once again integrate by parts over $B_{1}(0)$:
\begin{align}
\ve^{\kappa}\partial_t \overline{vc^{\ve}}+\ve\nabla_x\cdot \overline{vv^Tc^{\ve}}+a(\overline{vc^{\ve}}-v_*^{\ve}\overline{c^{\ve}})=
 \frac{n}{n+1}\E[q]\overline{c^{\ve}}-\overline{vc^{\ve}}.\label{M1e}
\end{align}
Rearranging and dividing \cref{M1e} by $\ve^{\kappa-1}$ leads to 
\begin{align}
 -(a+1)\ve^{1-\kappa}\overline{vc^{\ve}}=&\ve^{2-\kappa} \nabla_x\cdot \overline{vv^Tc^{\ve}}-\ve^{1-\kappa}\left(av_*^{\ve}+\frac{n}{n+1}\E[q]\right)\overline{c^{\ve}}+\ve\partial_t \overline{vc^{\ve}}.\label{sc2}
\end{align}
\noindent
Next, we apply $(\nabla_x\cdot)$ to both sides of \cref{sc2} and plug the expression on the right-hand side into \cref{mom0}. In order to eliminate the resulting term with the mixed derivative $(\nabla_x\cdot)\partial_t$ we apply ${\varepsilon^{\kappa}}\partial_t$ to both sides of \cref{mom0}. Thus we arrive at the following differential equation for the moments of zero and second order:
\begin{align}
 \ve^{\kappa}\partial_{t^2} \overline{c^{\ve}}+(a+1)\partial_t \overline{c^{\ve}}
 =&\ve^{2-\kappa}\nabla_x\nabla_x^T:\overline{vv^T c^{\ve}}-\ve^{1-\kappa}\nabla_x\cdot\left(\left(av_*^{\ve}+\frac{n}{n+1}\E[q]\right)\overline{c^{\ve}}\right).\label{mom012}
\end{align}
\begin{Remark}{ At first glance, equation \cref{mom012} is a damped wave equation with a transport term. Yet this is only the case if $\kappa=1$ and  $c^{\ve}$ is independent of $v$, so that the first term on the right-hand side becomes $\ve\Delta_x \overline{c^{\ve}}$. }
     \end{Remark}

\subsection{Zero order approximation}
Passing formally to the limit as $\ve\rightarrow0$ in \cref{transpce} and using \cref{vstare} we obtain that $$c^{0}:=\underset{\ve\rightarrow0}{\lim} \, c^{\ve}$$ satisfies the {equation}
\begin{align}
 -a\nabla_v\cdot(vc^{0})=nq\overline{c^{0}}-c^{0},\label{eqmu}
\end{align}
which can be solved explicitly:
\begin{Proposition}\label{prop-erste}
For any fixed $\overline{c^{0}}$ there exists a unique  solution to \cref{eqmu} with $c^{0}=0$ for $|v|=1$:
\begin{align}
 c^{0}=\overline{c^{0}}q\Cr{xi1},\label{solc0}
 \end{align}
 where
 {\begin{align}
 \Cl[xi]{xi1}(v)=\begin{cases}
 \frac{n}{na-1}\left(|v|^{-n+\frac{1}{a}}-1\right)&\text{for }a\neq \frac{1}{n},\\
 -\frac{n}{a}\ln|v|&\text{for }a=\frac{1}{n}.\end{cases}\label{xi1l}
\end{align}}
\end{Proposition}
\begin{proof}
We use the method of characteristics which leads to the following ODE system:
\begin{subequations}
\begin{alignat}{3}
 &\partial_{\sigma}v=-av,&&\qquad|v(0)|=1,\label{ODEv}\\
 &\partial_{\sigma}c^{0}=nq\overline{c^{0}}+(na-1)c^{0},&&\qquad c^{0}(0)=0.\label{ODEc0}
\end{alignat}
\end{subequations}
The ODE \cref{ODEv} leads to
\begin{align}
 v({\sigma})=e^{-a{\sigma}}v(0),\label{solvs}
\end{align}
so that
\begin{align}
 |v({\sigma})|=e^{-a{\sigma}}.\label{solvsn}
\end{align}
{Let $a\neq \frac{1}{n}$.} 
Solving the ODE \cref{ODEc0}, we obtain using \cref{qsign} and \cref{solvs} {that
\begin{align}
 c^{0}({\sigma})=&\overline{c^{0}}\int_0^{\sigma}\,ne^{(na-1)({\sigma}-\tau)}q(\hat v)(\tau)d\tau\nonumber\\
 =&\overline{c^{0}}q(\hat v)({\sigma})\int_0^{\sigma}\,ne^{(na-1)({\sigma}-\tau)}d\tau \nonumber\\
 =&\overline{c^{0}}q(\hat v)({\sigma})\frac{n}{na-1}\left(e^{(na-1){\sigma}}-1\right),\label{solvsc0}
\end{align}
}which together with \cref{solvsn} gives \cref{solc0}
{for $a\neq \frac{1}{n}$. Passing to the limit as  $a\rightarrow\frac{1}{n}$ gives the formula for $a=\frac{1}{n}$.}
\end{proof}

\noindent

\noindent
Next, we multiply \cref{eqmu} by $v$ and $vv^T$, respectively, and integrate by parts over $B_{1}(0)$ in order to find the  moments of order one:
\begin{align}
 &a\overline{vc^{0}}=\frac{n}{n+1}\E[q]\overline{c^{0}}-\overline{vc^{0}}\nonumber\\
 \Leftrightarrow\qquad & \overline{vc^{0}}=\frac{1}{a+1}\frac{n}{n+1}\E[q]\overline{c^{0}},\label{M1}
\end{align}
and two, respectively:
\begin{align}
 &2a\overline{vv^Tc^{0}}=n\overline{vv^Tq}\,\overline{c^{0}}-\overline{vv^Tc^{0}}\nonumber\\
 \Leftrightarrow\qquad &
 \overline{vv^Tc^{0}}=\frac{1}{2a+1}\frac{n}{n+2}\D[q]\overline{c^{0}}.\label{M2}
\end{align}
Passing formally to the limit in \cref{mom012}, using \cref{vstare} and \cref{M2} we arrive at  a  drift-diffusion equation
\begin{align}
(a+1)\partial_t\overline{c^{0}}=\frac{1}{2a+1}\frac{n}{n+2}\nabla_x\nabla_x^T: \left(\D[q]\overline{c^{0}}\right)-a\nabla_x\cdot(\overline{c^{0}} \F\nabla_x Q)\qquad\text{if }\kappa=2\text{ and }\E[q]\equiv0.\label{PL}
\end{align}
Passing formally to the limit in \cref{mom0} and using \cref{M1} we arrive at a drift equation 
\begin{align}
 (a+1)\partial_t \overline{c^{0}}+\frac{n}{n+1}\nabla_x\cdot(\overline{c^{0}}\E[q])=0\qquad\text{if }\kappa=1.\label{HL}
\end{align}
Several remarks on the resulting equations are in order.
\begin{Remark}
 In the parabolic scaling case ($\kappa =2$) the first term on the right hand side {of equation \cref{PL}} represents 
the myopic diffusion which
, if the scaling constant depending on $n$ and $a$ is neglected, can be decomposed as follows:
 \begin{align}
 \nabla_x\nabla_x^T: \left(\D[q]\overline{c^{0}}\right)=\nabla_x \cdot \left(\D[q]\nabla _x\overline{c^{0}}\right)+\nabla_x \cdot\left(\overline{c^{0}}\nabla _x\cdot \D[q]\right),\nonumber
 \end{align}
 where the first summand  is the regular anisotropic diffusion in divergence form with diffusion coefficient  $\D[q]$, and the second summand describes  cell transport  with velocity $\nabla _x\cdot \D[q]$.  {Due to \eqref{VD} we have that}  $\D[q]$ coincides with the variance-covariance matrix $\V[q]$. Both effects are thus dependent on the orientation distribution of tissue fibers. The second term on the right-hand side of \cref{PL} describes 
 the taxis towards the gradient of the macroscopic quantity $Q$ (e.g., haptotaxis if $Q$ represents macroscopic tissue density, chemotaxis if it denotes the concentration of some chemoattractant). When a hyperbolic scaling ($\kappa=1$) is performed, however, the dynamics is dominated by transport in the mean fiber direction - at least at 
 leading order.  In the next \cref{SecCorr} we will derive first order corrections for {\it both} scalings, which allows for  more accurate approximations.
 \end{Remark}
 \begin{Remark}[The role of $a$]
Both scaling limits involve a constant $a$, 
which is a scaling parameter for the acceleration. 
Specifically, $1/a$ can be seen to be analogous to the parameter employed in \cite{CEKNSSW} to characterise single cell velocity dynamics: it should be a quantity of the order $\frac{1}{\gamma}\epsilon ^{-\gamma}$, for some $\gamma >0$ representing a constant related to smaller scales, e.g. microtubule extension zones that are responsible for the subcellular level exchange of cells with their environment.  For further details we refer to \cite{CEKNSSW}. 
For both types of  scaling,  sending $a$ to zero leads to the standard diffusion and drift equations, 
which were previously derived for cell movement without transport with respect to velocity included in the KTE, see e.g.  \cite{HillenPainter}.  Conversely, for large $a$ the role of terms depending on the (mesoscopic) fiber orientation  becomes negligible. In the parabolic limit the taxis with respect to the macroscopic quantity $Q$ then fully dominates the space-time evolution. In the hyperbolic scaling case no such additional effect is present, so that the macroscopic cell density remains nearly constant over time.
%
\end{Remark}
\begin{Remark}[Directed/undirected fibers]\label{DirUndir}
As in previous studies, see, e.g. \cite{Hillen2006}, the parabolic scaling can only be performed under the assumption $\mathbb E[q]\equiv 0$. This occurs,  but not exclusively, in the case  where the tissue fibers are  undirected.  By this we mean (as, e.g. in \cite{Hillen2006}) that the fibers are symmetrical all along their axes, i.e. there is no 'up' and 'down' on such fibers, which
translates into symmetry of the orientational distribution: $$q(x,\theta )=q(x,-\theta)\qquad\text{for all }\theta \in S_1(0).$$
On the other hand, the hyperbolic scaling required no such assumption. 

As described in \cite{HillenPainter}, the choice of an appropriate scaling can be made based on measurements of reference values (such as speeds, turning rates, etc.) in a specific application. As far as 
brain tissue is concerned, it is still not clearly established whether it is directed or not, however recent mathematical modelling and simulations of typical glioblastoma patterns suggests that it might be undirected \cite{Kumar20}.
\end{Remark}
\subsection{First order correction}\label{SecCorr}
The above formal passage to the limit for $\varepsilon \to 0$ has led to macroscopic PDEs only containing leading order terms. To obtain first order corrections, hence  enhanced approximations, we start by introducing 
\begin{align*}
  &c^{0}_1:={\underset{\ve\rightarrow0}{\lim} \,\partial_{\ve}c^{\ve}},\\
  &c_{01}^{\ve}:=c^{0}+\ve c^{0}_1.
\end{align*}
Differentiating \cref{transpce} and  \cref{mom012}  with respect to $\ve${, letting $\ve\rightarrow0$,} and using \cref{vstare}  we obtain{: from \cref{transpce} that}
\begin{align}
 -a\nabla_v\cdot(vc_1^0)=nq\overline{c_1^0}-c_1^0-\left(\delta_{1\kappa}\partial_t c^{0}+\nabla_x\cdot (vc^{0})+a\F\nabla_x Q\cdot \nabla_v c^{0}\right),\label{eqc11}
\end{align}
{where $\delta_{1\kappa}$ denotes the Kronecker delta, and from   \cref{mom012} that} 
\begin{align}
 (a+1)\partial_t \overline{c^0_1}
 =&\nabla_x\nabla_x^T:\overline{vv^T c^0_1}-a\nabla_x\cdot\left(\overline{c^{0}_1}\F\nabla_x Q\right)+a\nabla_x\cdot\left(\overline{c^{0}}\F\nabla_x Q|\nabla_x Q|\right)\qquad\text{if }\kappa=2\text{ and }\E[q]\equiv0,\label{demom012P}
\end{align}
\begin{align}
 \partial_{t^2} \overline{c^0}+(a+1)\partial_t \overline{c^0_1}
 =&\nabla_x\nabla_x^T:\overline{vv^T c^0}-\nabla_x\cdot\left(\frac{n}{n+1}\E[q]\overline{c^0_1}+a\overline{c^0}\F\nabla_x Q\right)\qquad\text{if }\kappa=1.\label{demom012H}
\end{align}
\begin{Proposition}\label{PropHypc01}
For any  $\overline{c_1^0}$ there exists a unique  solution to \cref{eqc11} with $c_1^0=0$ for $|v|=1$:
\begin{align}
 c_1^0=\overline{c^{0}_1}q\Cr{xi1}-\delta_{1\kappa}\partial_t\overline{c^{0}}q\Cl[xi]{xi3}-\nabla_x\cdot\left(vq\overline{c^{0}}\right) \Cl[xi]{xi4}-a\overline{c^{0}}\F\nabla_x Q\cdot\left( \nabla_v q\Cl[xi]{xi5}+vq\Cl[xi]{xi6}\right),\label{solc10}
\end{align}
{where
\begin{align}
    &\Cr{xi3}(v)=\begin{cases}
 \frac{1}{a}\frac{n}{na-1}|v|^{-n+\frac{1}{a}}\left(-\ln|v|+\frac{1}{n-\frac{1}{a}}\left(|v|^{n-\frac{1}{a}}-1\right)\right)&\text{for }a\neq \frac{1}{n},\\
 \frac{n^3}{2}\ln^2|v|&\text{for }a=\frac{1}{n},\end{cases}
    \label{xxi3}\\
    &\Cr{xi4}(v)=\begin{cases}
  \frac{1}{a}\frac{n}{na-1}|v|^{-n+\frac{1}{a}-1}\left(1-|v|+\frac{1}{n-\frac{1}{a}+1}\left(|v|^{n-\frac{1}{a}+1}-1\right)\right)&\text{for }a\not\in\left\{\frac{1}{n+1},\frac{1}{n}\right\},\\
 n(n+1)^2(-\ln|v|-1+|v|)&\text{for }a=\frac{1}{n+1},\\
 n^3|v|^{-1}(1-|v|+|v|\ln|v|)&\text{for }a=\frac{1}{n},\end{cases}
   \\
    &\Cr{xi5}(v)=\begin{cases}
 \frac{1}{a}\frac{n}{na-1}|v|^{-n+\frac{1}{a}+1}\left(|v|^{-1}-1+\frac{1}{n-\frac{1}{a}-1}\left(|v|^{n-\frac{1}{a}-1}-1\right)\right)&\text{for }a\not\in\left\{ \frac{1}{n-1},\frac{1}{n}\right\},\\
 n(n-1)^2|v|^{-1}(1-|v|+|v|\ln|v|)&\text{for }a=\frac{1}{n-1}\text{ and }n>1,\\
 n^3(-\ln|v|-1+|v|)&\text{for }a=\frac{1}{n},\end{cases}
    \\
    &\Cr{xi6}(v)=
 \frac{n}{a^2}\left(|v|^{-1}-1\right).
    \label{xxi6}
 \end{align}}
\end{Proposition}
\begin{Remark} Notice that the first order correction $c_1^0$ and leading term {$c^0$} 
given by  \cref{solc10} and  \cref{solc0}, respectively, are clearly not mutually orthogonal with respect to the scalar product of the weighted $L^2$-space $L^2(V;\frac{dv}{nq})$.
\end{Remark}

\begin{proof}(of \cref{PropHypc01}) 
Observe first that  \cref{qsign} implies 
\begin{align}
 {\nabla_v q(\hat v)
 =({\mathbb I_n}-{\hat v \hat v^T})\nabla_vq(\hat v)|v|^{-1}.}\label{nqsign}
\end{align}
Hence
\begin{align}
 \nabla_v(q(\hat v)\Cr{xi1}(|v|))=\left({\mathbb I_n}-{\hat v \hat v^T}\right)\nabla_vq(\hat v)|v|^{-1}\Cr{xi1}(|v|)+q(\hat v)\hat v\Cr{xi1}'(|v|).
 \label{nabc0}
\end{align}
Using the method of characteristics and  \cref{nabc0,solvs,solc0} we obtain that
 \begin{align}
  c^{0}_1(\sigma)=&\int_0^\sigma\,e^{(na-1)(\sigma-\tau)}\left(nq\overline{c_1^0}-\left(\delta_{1\kappa}\partial_t c^{0}+\nabla_x\cdot (vc^{0})+a\F\nabla_x Q\cdot \nabla_vc^{0}\right)\right)(\tau)\,d\tau\nonumber\\
  =&\overline{c^{0}_1}q\Cr{xi1}-\int_0^\sigma\,e^{(na-1)(\sigma-\tau)}\left(\delta_{1\kappa}\partial_t c^{0}+\nabla_x\cdot (vc^{0})+a\F\nabla_x Q\cdot \nabla_vc^{0}\right)\,d\tau\nonumber\\
  =&\overline{c^{0}_1}q\Cr{xi1}-\delta_{1\kappa}\partial_t\overline{c^{0}}q\Cr{xi3}-\nabla_x\cdot\left(vq\overline{c^{0}}\right) \Cr{xi4}-a\overline{c^{0}}\F\nabla_x Q\cdot\left( \nabla_v q\Cr{xi5}+vq\Cr{xi6}\right),
 \end{align}
 {where 
 \begin{align}
    &\Cr{xi3}(v)(\sigma)=e^{(na-1)\sigma}\int_0^\sigma\,e^{-(na-1)\tau}\Cr{xi1}(|v|)(\tau)\,d\tau,\label{xxi3_1}\\
    &\Cr{xi4}(v)(\sigma)=|v|^{-1}(\sigma)e^{(na-1)\sigma}\int_0^\sigma\,e^{-(na-1)\tau}|v|\Cr{xi1}(|v|)(\tau)\,d\tau,\\
    &\Cr{xi5}(v)(\sigma)=|v|(\sigma)e^{(na-1)\sigma}\int_0^\sigma\,e^{-(na-1)\tau}|v|^{-1}\Cr{xi1}(|v|)(\tau)\,d\tau,\\
    &\Cr{xi6}(v)(\sigma)=e^{(na-1)\sigma}\int_0^\sigma\,e^{-(na-1)\tau}\Cr{xi1}'(|v|)(\tau)\,d\tau.\label{xxi6_1}
 \end{align}
 Since $e^{-(na-1)\tau}=|v|^{n-\frac{1}{a}}(\tau)$ and  $e^{-(na-1)\tau}\,d\tau=-\frac{1}{a}|v|^{n-\frac{1}{a}-1}\,d|v|$ due to \cref{solvsn}, we can rewrite \cref{xxi3_1}-\cref{xxi6_1} as follows:
 \begin{align}
    &\Cr{xi3}(v)=\frac{1}{a}|v|^{-n+\frac{1}{a}}\int^1_{|v|}\,s^{n-\frac{1}{a}-1}\Cr{xi1}(s)\,ds,\label{xxi3_2}\\
    &\Cr{xi4}(v)=\frac{1}{a}|v|^{-n+\frac{1}{a}-1}\int^1_{|v|}\,s^{n-\frac{1}{a}}\Cr{xi1}(s)\,ds,\\
    &\Cr{xi5}(v)=\frac{1}{a}|v|^{-n+\frac{1}{a}+1}\int^1_{|v|}\,s^{n-\frac{1}{a}-2}\Cr{xi1}(s)\,ds,\\
    &\Cr{xi6}(v)=\frac{1}{a}|v|^{-n+\frac{1}{a}}\int^1_{|v|}\,s^{n-\frac{1}{a}-1}\Cr{xi1}'(s)\,ds.\label{xxi6_2}
 \end{align}
 Let $a\neq\frac{1}{n}$. Plugging \cref{xi1l} into  \cref{xxi3_2}-\cref{xxi6_2}, we obtain
 \begin{align}
    &\Cr{xi3}(v)=\frac{1}{a}\frac{n}{na-1}|v|^{-n+\frac{1}{a}}\int^1_{|v|}\,s^{n-\frac{1}{a}-1}\left(s^{-n+\frac{1}{a}}-1\right)\,ds,\label{xxi3_3}\\
    &\Cr{xi4}(v)=\frac{1}{a}\frac{n}{na-1}|v|^{-n+\frac{1}{a}-1}\int^1_{|v|}\,s^{n-\frac{1}{a}}\left(s^{-n+\frac{1}{a}}-1\right)\,ds,\\
    &\Cr{xi5}(v)=\frac{1}{a}\frac{n}{na-1}|v|^{-n+\frac{1}{a}+1}\int^1_{|v|}\,s^{n-\frac{1}{a}-2}\left(s^{-n+\frac{1}{a}}-1\right)\,ds,\\
    &\Cr{xi6}(v)=-\frac{n}{a^2}|v|^{-n+\frac{1}{a}}\int^1_{|v|}\,s^{-2}\,ds.\label{xxi6_3}
 \end{align}
 Finally, computing the integrals in \cref{xxi3_2}-\cref{xxi6_2}, we arrive at \cref{xxi3}-\cref{xxi6} for $a\neq \frac{1}{n}$. Passing to the limit as  $a\rightarrow\frac{1}{n}$ gives the formulas for $a=\frac{1}{n}$.
 }
\end{proof}
\noindent
Next, we multiply \cref{eqc11} by $v$ and $vv^T$, respectively, integrate by parts over $B_{1}(0)$, and use \cref{M1,M2} in order to find the  moments of order one:
\begin{align}
 &a\overline{vc^{0}_1}=\frac{n}{n+1}\E[q]\overline{c^{0}_1}-\overline{vc^{0}_1}-\left(\delta_{1\kappa}\partial_t \overline{vc^{0}}+\nabla_x\cdot \overline{vv^Tc^{0}}-a\overline{c^{0}}{\F}\nabla_x Q \right)\nonumber\\
 \Leftrightarrow\quad &(a+1)\overline{vc^{0}_1}=\frac{n}{n+1}\E[q]\overline{c^{0}_1}-\left(\delta_{1\kappa}\frac{1}{a+1}\frac{n}{n+1}\E[q]\partial_t \overline{c^{0}}+\frac{1}{2a+1}\frac{n}{n+2}\nabla_x\cdot \left(\D[q]\overline{c^{0}}\right)-a\overline{c^{0}}{\F}\nabla_x Q \right),\label{M1cor}
\end{align}
and two,  respectively:
\begin{align}
 &2a\overline{vv^Tc^{0}_1}=n\overline{vv^Tq}\,\overline{c^{0}_1}-\overline{vv^Tc^{0}_1}-\left(\delta_{1\kappa}\partial_t \overline{vv^Tc^{0}}+\overline{vv^T(v\cdot\nabla_x c^{0})}{-}a\left({\F}\nabla_x Q\overline{vc^{0}}^T+\overline{vc^{0}} ({\F}\nabla_x Q)^T\right)\right)\nonumber\\
 \Leftrightarrow\quad &
 (2a+1)\overline{vv^Tc^{0}_1}=\frac{n}{n+2}\D[q]\overline{c^{0}_1}
 -\Bigg(\frac{\delta_{1\kappa}}{2a+1}\frac{n}{n+2}\D[q]\partial_t \overline{c^{0}}{+\Cl{C10}\nabla _x\cdot \left (\T[q]\ \overline{c^{0}}\right)}\nonumber\\ 	
 &\qquad\qquad\qquad\quad+\frac{a}{a+1}\frac{n}{n+1}\overline{c^{0}}\left({\F}\nabla_x Q\ \E[q]^T+\E[q]({\F}\nabla_x Q)^T\right)\Bigg),\label{M2cor}
\end{align}
where $C_1$ is a constant depending on $a$ and $n$.

\subsubsection{Parabolic scaling}
Let $\kappa=2$ and $\E[q]\equiv 0$, $\T[q]\equiv 0$. 
Then \cref{M2cor} simplifies to
\begin{align}
 &\overline{vv^Tc^{0}_1}=\frac{1}{2a+1}\frac{n}{n+2}\D[q]\overline{c^{0}_1}.\label{M2cork2}
\end{align}
Plugging \cref{M2cork2} into \cref{demom012P} we arrive at an equation for $\overline{c^{0}_1}$:
\begin{align}
 (a+1)\partial_t \overline{c^{0}_1}
 =&\frac{1}{2a+1}\frac{n}{n+2}\nabla_x\nabla_x^T: \left(\D[q]\overline{c^{0}_1}\right)-a\nabla_x\cdot\left(\overline{c^{0}_1}\F\nabla_x Q\right)+a\nabla_x\cdot\left(\overline{c^{0}}\F\nabla_x Q|\nabla_x Q|\right).\label{dif2}
\end{align}
Combining \cref{PL,dif2} and using \cref{vstare}, we obtain for $\overline{c_{01}^{\ve}}$ the equation
\begin{align}
 (a+1)\partial_t \overline{c_{01}^{\ve}}
 =\frac{1}{2a+1}\frac{n}{n+2}\nabla_x\nabla_x^T: \left(\D[q]\overline{c_{01}^{\ve}}\right)-a\nabla_x\cdot\left(\overline{c_{01}^{\ve}}\F\frac{\nabla_x Q}{1+\ve|\nabla_x Q|}\right)+O\left(\ve^2\right)\quad&\text{if } \kappa=2\text{ and}\nonumber\\
 &\E[q]\equiv0,\   \T[q]\equiv0.\label{dif3}
\end{align}\begin{Remark}
 A first order correction is often neglected in parabolic scaling. Here it allows to get a description of the taxis with respect to the macroscopic quantity $Q$ that is more accurate than in \cref{PL}. 
It turns out to be closer to a flux-limited taxis. The  myopic diffusion remains unchanged.
\end{Remark}
\cgr{\begin{Remark}
      Both conditions $\E[q]\equiv0$ and $   \T[q]\equiv0$  are automatically satisfied  if { the} tissue fibers are undirected (see \cref{DirUndir}). 
     \end{Remark}
}
\subsubsection{Hyperbolic scaling}
Let $\kappa=1$. Plugging \cref{M2} into  \cref{demom012H}  we arrive at the equation
\begin{align}
 (a+1)\partial_t \overline{c^{0}_1}
 +\frac{n}{n+1}\nabla_x\cdot\left(\E[q]\overline{c^{0}_1}\right)=&-\partial_{t^2} \overline{c^{0}}+\frac{1}{2a+1}\frac{n}{n+2}\nabla_x\nabla_x^T:\left(\D[q]\overline{c^{0}}\right)-a\nabla_x\cdot\left(\overline{c^{0}}\F\nabla_xQ\right).\label{dif2H}
\end{align}
Utilising \cref{HL} twice we compute 
\begin{align}
 -\partial_{t^2} \overline{c^{0}}=&\frac{1}{a+1}\frac{n}{n+1}\nabla_x\cdot\left(\partial_t\overline{c^{0}}\E[q]\right)\nonumber\\
 =&-\frac{1}{(a+1)^2}\frac{n^2}{(n+1)^2}\nabla_x\cdot\left(\E[q]\nabla_x\cdot\left(\overline{c^{0}}\E[q]\right)\right).\label{difft2}
\end{align}
Plugging \cref{difft2} into \cref{dif2H} we obtain that
\begin{align}
 &(a+1)\partial_t \overline{c^{0}_1}
 +\frac{n}{n+1}\nabla_x\cdot\left(\E[q]\overline{c^{0}_1}\right)\nonumber\\
 =&\frac{1}{2a+1}\frac{n}{n+2}\nabla_x\nabla_x^T:\left(\D[q]\overline{c^{0}}\right)-\frac{1}{(a+1)^2}\frac{n^2}{(n+1)^2}\nabla_x\cdot\left(\E[q]\nabla_x\cdot\left(\overline{c^{0}}\E[q]\right)\right)-a\nabla_x\cdot\left(\overline{c^{0}}\F\nabla_xQ\right).\label{dif2H2}
\end{align}
Consequently,  in virtue of \cref{HL} we obtain that  $\overline{c_{01}^{\ve}}$ satisfies the equation
\begin{align}
 &(a+1)\partial_t \overline{c_{01}^{\ve}}
 +\frac{n}{n+1}\nabla_x\cdot\left(\E[q]\overline{c_{01}^{\ve}}\right)\nonumber\\
 =& {\ve}\left(\frac{1}{2a+1}\frac{n}{n+2}\nabla_x\nabla_x^T:\left(\D[q]\overline{c_{01}^{\ve}}\right)-\frac{1}{(a+1)^2}\frac{n^2}{(n+1)^2}\nabla_x\cdot\left(\E[q]\nabla_x\cdot\left(\overline{c_{01}^{\ve}}\E[q]\right)\right)-a\nabla_x\cdot\left(\overline{c_{01}^{\ve}}\F\nabla_xQ\right)\right)\nonumber\\
 &+O\left(\ve^2\right).\label{dif4H}
\end{align}
\begin{Remark}
Typically for  the hyperbolic case,  the first order correction includes two terms which depend on the mesoscopic fiber orientation distribution: a myopic diffusion (the same as in  \cref{dif3}) and yet another term contributing to cell diffusion as well as to transport. The resulting diffusion tensor is  (up to multiplication by a constant) a linear combination of two nonnegative definite matrices:
 \begin{align*}
  \mathbb V[q]+\left (1-\frac{(2a+1)n(n+2)}{(a+1)^2(n+1)^2}\right )\E[q]\otimes \E[q].
 \end{align*}
 Since the coefficient before the second matrix is obviously a positive number for any $a>0$ and $n\in\N$, the  diffusion tensor is nonnegative definite. 
 \\\indent 
 The final term in the second line of \cref{dif4H} describes taxis with respect to the macroscopic quantity $Q$. As is the case for the  parabolic limit equation \cref{PL}, the resulting taxis is not flux-limited. Even higher order approximations  are necessary in order to capture this effect properly.

 Overall, the first order correction effects a (small) deviation from the drift towards the average orientation of tissue fibers.
\end{Remark}

\begin{Remark}[Approximation order]\ \\[-2.5ex] \begin{enumerate} \item Equation  \cref{dif4H} shows that the error of the first order approximation is  $O(\ve^2)$. This confirms the surmise stated in \cite{HillenPainter} in connection with the hyperbolic scaling  performed for a closely related problem. There, however, a different approach which is based on the  Chapman-Enskog expansion was used in order to approximate the mesoscopic cell density.
\item
Differentiating \cref{transpce} and \cref{mom012}  with respect to $\ve$ at $\ve=0$ several times and performing the necessary calculations in the same manner as above, one obtains equations for corrections of higher order.
\end{enumerate}

\end{Remark}



\subsection{A mesoscopic first order approximation}\label{SecMesoAppr}
One known issue with the above approach relying on {the} first order approximation is that $c_{01}^{\ve}$  {is not necessarily nonnegative even if $\overline{c^{0}_1}$ is nowhere negative.} To fix this drawback one could consider instead an approximation $\widetilde{c_{01}^{\ve}}$ which {vanishes for $v\in\partial B_{1}$ and} solves the following equation:
\begin{subequations}\label{TransTilde}
\begin{align}
 &\ve^{\kappa}\partial_t \widetilde{c_{01}^{\ve}}+\ve \nabla_x\cdot \left(v\widetilde{c_{01}^{\ve}}\right)-a\nabla_v\cdot\left((v-{v_*^{\ve}})\widetilde{c_{01}^{\ve}}\right)+\widetilde{c_{01}^{\ve}}=nq\overline{c_{01}^{\ve}}=nq{(}\overline{c^{0}+\ve c^{0}_1}{)},\label{transpce1}\\
 &{\widetilde{c_{01}^{\ve}}|_{t=0}=c_0.}\label{ctildeini}
\end{align}
\end{subequations}
{Equation \cref{transpce1} is obtained from the original KTE \cref{transpce} upon replacing the integral term by its first order approximation.} 
Our next Proposition shows that the two main properties  any reasonable approximation of a mesoscopic density should have are satisfied: both positivity and the total mass are preserved under \cref{transpce1}.
\begin{Proposition}
      Let  $\widetilde{c_{01}^{\ve}}$ be a solution to \cref{TransTilde}. Suppose that   $\overline{c_{01}^{\ve}}$ and $c_0$ are nonnegative. Then $\widetilde{c_{01}^{\ve}}$ is also nonnegative and satisfies
      \begin{align}
       \int_{\R^n}\overline{\widetilde{c_{01}^{\ve}}}\,dx\equiv\int_{\R^n}\overline{c_0}\,dx,\label{TildeMas}
      \end{align}
      i.e. the total mass is preserved.      
     \end{Proposition}
\begin{proof}
To begin with, we apply the method of characteristics to the transport equation \cref{transpce1} which  leads to the ODE system
\begin{subequations}
\begin{align}
  \ve^{\kappa-1}\partial_{\sigma}x=&v,\\
  \ve^{\kappa}\partial_{\sigma}v=&-a(v-{v_*^{\ve}}(x)),\\
 \ve^{\kappa}\partial_{\sigma}\widetilde{c_{01}^{\ve}} -(na-1)\widetilde{c_{01}^{\ve}}=&nq\overline{c_{01}^{\ve}}.\label{ctildeODE}
\end{align}
\end{subequations}
It is obvious from \cref{ctildeODE} that if  $\overline{c_{01}^{\ve}}$ and $c_0$ are nonnegative, then  $\widetilde{c_{01}^{\ve}}$ is nonnegative as well. It remains to prove \cref{TildeMas}. Recall that the full  mesoscopic density $c^{\ve}$ solves the measure conserving KTE \cref{transpce} (this follows by integration of  \cref{mom0} by parts over $\R^n$ with respect to $x$), so that
\begin{align}
 \int_{\R^n}\overline{c^{\ve}}\,dx\equiv&\int_{\R^n}\overline{c_{0}}\,dx.\nonumber
\end{align}
This and the definition of $c_{01}^{\ve}$ entail 
\begin{align}
 \int_{\R^n}\overline{c_{01}^{\ve}}\,dx\equiv&\int_{\R^n}\overline{c_{0}}\,dx.\label{c01totm}
\end{align}
Further, we integrate  \cref{transpce1} with respect to $x$ and $v$ over the whole space $\R^n\times \overline{B}_{1}(0)$. Using  partial integration and \cref{c01totm} we obtain an ODE for the total mass:
\begin{align}
 \frac{d}{dt}\int_{\R^n}\overline{\widetilde{c_{01}^{\ve}}}\,dx+\int_{\R^n}\overline{\widetilde{c_{01}^{\ve}}}\,dx=&\int_{\R^n}\overline{c_{01}^{\ve}}\,dx,\nonumber\\
 =&\int_{\R^n}\overline{c_{0}}\,dx.\label{ctildeODEMas}
\end{align}
Finally, we integrate the initial condition \cref{ctildeini} over the whole space and obtain
\begin{align}
 \int_{\R^n}\overline{\widetilde{c_{01}^{\ve}}}\,dx|_{t=0}=&\int_{\R^n}\overline{c_0}.\label{ctildeiniMas}
\end{align}
Solving \cref{ctildeODEMas}-\cref{ctildeiniMas} yields \cref{TildeMas}.
\end{proof}
\begin{Remark}
  In the previous Subsection we have  obtained {DTE}s for $\overline{c_{01}^{\ve}}$ for both parabolic and hyperbolic scaling cases.  Each of those macroscopic equations for $\overline{c_{01}^{\ve}}$ can be solved numerically,  so that \cref{transpce1} can be regarded  as a linear transport equation which needs to be solved in order to determine $\widetilde{c_{01}^{\ve}}$.  
    This seems to be a useful alternative to dealing  directly  with the KTE \cref{transpce1}, since  the numerical handling of such equations is known to be more expensive. 
\end{Remark}
\section{Rigorous limit procedures}\label{SecRig}
\subsection{Functional spaces}\label{SecFSp}
 {We begin with some basic notation. 
Let $O$ be a domain or a smooth manifold. We denote by $C^k_b(\overline{O})$ the space of real-valued functions which are  continuous and bounded together with their derivatives up to order $k$. For $k=0$ we write $C_b(\overline{O})$. If $\overline{O}$ is compact, we suppress the index $b$. For a set $A\subset \overline{O}$ we denote by $C^k_c(A)$ the subset  of $C^k_b(\overline{O})$ which consists of functions which are compactly supported in $A$.
Similarly, $C_b(\overline{O};X)$ for $X$ a Banach space denotes the space of continuous and bounded maps between $\overline{O}$ and $X$.}

Now we introduce various spaces of measures. We denote by ${\cal M}(\overline{O})$  the Banach space of finite signed Radon measures in $\overline{O}$, while   ${\cal M}_+(\overline{O})$ stands for the closed subset of  positive Radon measures. {As usual, we use the total variation as norm on ${\cal M}(\overline{O})$.}

\noindent
The product $\varphi\mu$ of a Radon measure $\mu\in {\cal M}(\overline{O})$ and a  function \cgr{$\varphi\in C_{b}(\overline{O})$} is understood in the usual sense: it is a Radon measure which satisfies
\begin{align}
 \int_{\cgr{\overline{O}} }\psi(v)\,d(\varphi\mu)(v)=\int_{\cgr{\overline{O}} }\psi\varphi(v)\,d\mu(v)\qquad \text{for all }\psi\in C_{\cgr c}(\overline{O}).\nonumber
\end{align}
We recall that due to the  Riesz-Markov-Kakutani representation theorem ${\cal M}(\overline{O})$ is  isometrically isomorphic to the continuous dual of the separable \cgr{normed} space $C_{\cgr c}(\overline{O})$.
Further, we extend the $\bar{~}$-operator which denotes the integration over $\overline{B}_{1}(0)$ to the case of measures: for 
$\mu\in {\cal M}(\overline{B}_{1}(0))$ let 
\begin{align}
 \overline{\mu}:=\int_{\overline{B}_{1}(0)}\,d\mu(v).
\end{align}
In order to simplify the notation, we identify a measure $q\in {\cal M}(S_1(0))$ with the measure $q\times (r^{n-1}\,dr)\in {\cal M}(\overline{B}_{1}(0))$.

\noindent
We also make use of the Banach space
\begin{align*}
L^{\infty}_{w\text{-}*}(\R^+;{\cal M}(\overline{O})):=&\left\{\mu:\R^+\rightarrow {\cal M}(\overline{O})\quad\text{is weak-}\!*\text{ measurable and}\right.\\ 
&\ \left. \left\|\mu\right\|_{L^{\infty}_{w\text{-}*}(\R^+;{\cal M}(\overline{O}))}:=\left\|\|\mu(\cdot)\|_{{\cal M}(\overline{O})}\right\|_{L^{\infty}(\R^+)}<\infty\right\}
\end{align*}
and its closed subset
\begin{align*}
&L^{\infty}_{w\text{-}*}(\R^+;{\cal M}_+(\overline{O})):=\left\{\mu:\R^+\rightarrow {\cal M}_+(\overline{O})\ |\ \mu\in L^{\infty}_{w\text{-}*}(\R^+;{\cal M}(\overline{O}))\right\}.
\end{align*}
Thereby we identify functions which coincide a.e. in $\R^+$. It is known \cite[sections 8.18.1-8.18.2]{Edwards} that
$L^{\infty}_{w\text{-}*}(\R^+;{\cal M}(\overline{O}))$ is isometrically isomorphic to the continuous dual of the Bochner space $L^1(\R^+,C_{\cgr c}(\overline{O})))$ via the duality paring 
\begin{align*}
 \left<\mu,\varphi\right>=\int_{\R^+}\int_{\overline{O}}\varphi(t)(x)\, (d\mu(t))(x)dt.
\end{align*}
Since $L^1(\R^+,C_{\cgr c}(\overline{O})))$ is separable, the Banach-Alaoglu theorem implies that balls in $L^{\infty}_{w\text{-}*}(\R^+;{\cal M}(\overline{O}))$ are weak-$*$ sequentially compact. 
\indent{Finally, we introduce the spaces
\begin{align*}
C_{w\text{-}*}(\R^+;{\cal M}(\overline{O})):=&\left\{\mu:\R^+\rightarrow {\cal M}(\overline{O})\quad\text{is weak-}\!*\text{ continuous and}\right.\\ 
&\ \left. \left\|\|\mu(\cdot)\|_{{\cal M}(\overline{O})}\right\|_{L^{\infty}(\R^+)}<\infty\right\}
\end{align*}
and its closed subset}
\begin{align*}
&{C_{w\text{-}*}(\R^+;{\cal M}_+(\overline{O})):=\left\{\mu:\R^+\rightarrow {\cal M}_+(\overline{O})\ |\ \mu\in C_{w\text{-}*}(\R^+;{\cal M}(\overline{O}))\right\}.}
\end{align*}
\subsection{{Main results and their proofs}}
From now on we assume the model parameters to satisfy the following conditions:
\begin{Assumptions}\label{AssqQ}~
 \begin{enumerate}
  \item $q\in C_b(\R^n;{\cal M}_+(S_1(0)))$ and 
 \begin{align}
\overline{q}\equiv \frac{1}{n}; \nonumber
\end{align}
\item {$Q\in C_b(\R^+_0\times\R^n)$, $\nabla_x Q\in C_b(\R^+_0\times\R^n;\R^n)$,  $\F\in C_b(\R^n;\R^{n\times n})$,}
\begin{align*}
      \|\|\F\|_{2}\|_{C_b(\R^{n\times n})}\leq1.
     \end{align*}
 \end{enumerate}
\end{Assumptions}
\begin{Remark}[Moments of $q$]
 Due the regularity assumption on $q$
the moments
\begin{align}
 &\E[q](x):=\int_{S_1(0)}v\,(dq(x))(v),\\
 &\D[q](x):=\int_{S_1(0)}vv^T\,(dq(x))(v)
\end{align}
are well-defined  and satisfy
\begin{align}
 \E[q]\in C_b(\R^n{;\R^n}),\ \D[q]\in C_b(\R^n{;\R^n\times\R^n}).
\end{align}
\end{Remark}

\noindent
Next, we define weak measure-valued solutions to the KTE \cref{transpce} and the limit equations.
\begin{Definition}[Weak solutions to the KTE]\label{WKKTE}
 Let {$\ve>0$ and let}  \cref{AssqQ} be satisfied and let  $c^{\ve}_0\in {\cal M}_+(\R^n\times \overline{B}_{1}(0))$ be some initial data.
 We call an element  $c^{\ve}\in L^{\infty}_{w\text{-}*}(\R^+;{\cal M}_+(\R^n\times \overline{B}_{1}(0)))$ a weak solution to \cref{transpce}
 if for all $\varphi\in C^1_{\cgr{c}}(\R^n)$, $\psi\in C^1(\overline{B}_{1}(0))$, and $\eta\in C_{{c}}^1[0,\infty)$ it holds that
 \begin{align}
 &{-}\ve^{\kappa}\eta(0)\int_{\R^n\times\overline{B}_{1}(0)}\varphi\psi\,dc^{\ve}_0-\int_0^{\infty}\int_{\R^n\times\overline{B}_{1}(0)}(\ve^{\kappa},\ve v,-a(v-{v_*^{\ve}}))\cdot\nabla_{(t,x,v)}(\eta\varphi\psi)\,dc^{\ve}(t)\,dt\nonumber\\
 =&\int_0^{\infty}\eta\left(n\int_{\R^n}\varphi\int_{\overline{B}_{1}(0)}\psi\,dq\,d\overline{c^{\ve}}(t)-\int_{\R^n\times\overline{B}_{1}(0)}\varphi\psi\,dc^{\ve}(t)\right)dt\label{mesoweaksol},
\end{align}
{and the total mass is preserved:
\begin{align}
 \int_{\R^n}\,d\overline{c^{\ve}}
 =&\int_{\R^n}\,d\overline{c^{\ve}_0}\qquad\text{a.e. in }(0,\infty).\label{conv}
\end{align} 
}
\end{Definition}
\begin{Remark}[Solvability]\label{RemSolv}
In this work we are mostly interested in upscaling and do not deal with solvability of KTEs. These equations have been studied by many authors, though mostly in the physical context. A sketch of the proof of the existence of solutions to {\cref{meso}}  is provided in  \cref{Append} for the reader's convenience. 
\end{Remark}

\begin{Definition}[Weak solutions to the parabolic limit]\label{WKPL}
 Let  {\cref{AssqQ}} be satisfied and let  $\overline{c^{0}_{0}}\in {\cal M}_+(\R^n)$ be some initial data. We call an element  $\overline{c^{0}}\in L^{\infty}_{w\text{-}*}(\R^+;{\cal M}_+(\R^n))$ a weak solution to \cref{PL}
 if for all $\varphi\in C^2_{\cgr{c}}(\R^n)$ and $\eta\in C_{{c}}^1[0,\infty)$ it holds that
 \begin{align}
 &{-}\eta(0)\int_{\R^n}\varphi\,d\overline{c^{0}_{0}}-\int_0^{\infty}\frac{d\eta}{dt}\int_{\R^n}\varphi\,d\overline{c^{0}}(t)\,dt\nonumber\\
 =&\int_0^{\infty}\eta \left(\int_{\R^n}\nabla_x\nabla_x^T\varphi: \frac{1}{(a+1)(2a+1)}\frac{n}{n+2}\D[q]\,d \overline{c^{0}}(t)+\int_{\R^n}\nabla_x\varphi\cdot\frac{a}{a+1}\F\nabla_x Q\,d\overline{c^{0}}(t) \right)dt.\label{PLsol}
\end{align}
\end{Definition}
\begin{Definition}[Weak solutions to the hyperbolic limit]\label{WKHL}
 Let  {\cref{AssqQ}} be satisfied and let  $\overline{c^{0}_{0}}\in {\cal M}_+(\R^n)$ be some initial data. We call an element  $\overline{c^{0}}\in L^{\infty}_{w\text{-}*}(\R^+;{\cal M}_+(\R^n))$ a weak solution to \cref{HL}
 if for all $\varphi\in C^1_{\cgr{c}}(\R^n)$ and $\eta\in C_{{c}}^1[0,\infty)$ it holds that
 \begin{align}
 &{-}\eta(0)\int_{\R^n}\varphi\,d\overline{c^{0}_{0}}-\int_0^{\infty}\left(\frac{d\eta}{dt}\int_{\R^n}\varphi\,d\overline{c^{0}}(t)+\eta\int_{\R^n}\nabla_x\varphi\cdot\frac{1}{a+1}\frac{n}{n+1}\E[q]\,d\overline{c^{0}}(t)\right)dt=0.\label{HLsol}
\end{align}
\end{Definition}
\begin{Remark}
 Due to the assumptions made on $q$, $Q$, $\F$, and $c^{\ve}_0${,} each of the integrals in {\cref{mesoweaksol}-\cref{HLsol}} is well-defined and finite. 
\end{Remark}

\noindent
Thus defined weak solutions satisfy in a weak sense the equations for the moments which were formally derived in {\cref{SecScaling}}:
 \begin{Proposition}[Weak moment equations]\label{PropWMOME}
 Let $c^{\ve}$ be a weak solution as in {\cref{WKKTE}}. Then equations \cref{mom0} and \cref{mom012} are satisfied, { respectively,} in the following weak sense: for all $\varphi\in C^1_{\cgr{c}}(\R^n)$ and $\eta\in C_{{c}}^1[0,\infty)$
 \begin{align}
 &{-}\eta(0)\int_{\R^n}\varphi\,d\overline{c^{\ve}_0}-\int_0^{\infty}\left(\frac{d\eta}{dt}\int_{\R^n}\varphi\,d\overline{c^{\ve}}(t)+\ve^{1-\kappa}\eta\int_{\R^n}\nabla_x\varphi\cdot\,d\overline{vc^{\ve}}(t)\right)\,dt=0{,}\label{mom01w}
\end{align}
and for all $\varphi\in C^2_{\cgr{c}}(\R^n)$ and $\eta\in C_{{c}}^2[0,\infty)$
 \begin{align}
 &{-}\ve^{\kappa}\frac{1}{a+1}\frac{d\eta}{dt}(0)\int_{\R^n}\varphi\,d\overline{c^{\ve}_0}-\ve^{\kappa}\frac{1}{a+1}\int_0^{\infty}\frac{d^2\eta}{dt^2}\int_{\R^n}\varphi\,d\overline{c^{\ve}}(t)\,dt{-}\eta(0)\int_{\R^n}\varphi\,d\overline{c^{\ve}_0}-\int_0^{\infty}\frac{d\eta}{dt}\int_{\R^n}\varphi\,d\overline{c^{\ve}}(t)\,dt\nonumber\\
 =&\int_0^{\infty}\eta\int_{\R^n}\nabla_x\nabla_x^T\varphi:\ve^{2-\kappa} \frac{1}{a+1}\,d \overline{vv^Tc^{\ve}}(t)dt\nonumber\\
 &+\int_0^{\infty}\eta\left(\int_{\R^n}\nabla_x\varphi\cdot\ve^{1-\kappa}\frac{a}{a+1}{v_*^{\ve}}\,d\overline{c^{\ve}}(t)+\int_{\R^n}\nabla_x\varphi\cdot\ve^{1-\kappa}\frac{1}{a+1}\frac{n}{n+1}\E[q]\, d\overline{c^{\ve}}(t)\right)dt\nonumber\\
 &-\eta(0)\int_{\R^n}\nabla_x\varphi\cdot\ve\frac{1}{a+1}\,d\overline{vc_0^{\ve}}.\label{mom012w}
\end{align}
\end{Proposition}
\begin{proof}
 Taking $\psi\equiv 1$ in \cref{mesoweaksol} directly implies \cref{mom01w}. Similarly, using  $\psi(v)=v$ and $\nabla_x \varphi$  in place of $\varphi$ for some $\varphi\in C^2_{\cgr{c}}(\R^n)$, we {obtain from \cref{mesoweaksol} that}
 \begin{align}
 &\int_0^{\infty}\ve^{1-\kappa}\eta\int_{\R^n}\nabla_x\varphi\cdot d\overline{vc^{\ve}}(t)dt\nonumber\\
 =&\int_0^{\infty}\eta\int_{\R^n}\nabla_x\nabla_x^T\varphi:\ve^{2-\kappa} \frac{1}{a+1}\,d \overline{vv^Tc^{\ve}}(t)dt\nonumber\\
 &+\int_0^{\infty}\eta\int_{\R^n}\nabla_x\varphi\cdot\left(\ve^{1-\kappa}\frac{a}{a+1}{v_*^{\ve}}\,d\overline{c^{\ve}}(t)+\ve^{1-\kappa}\frac{1}{a+1}\frac{n}{n+1}\,\E[q] d\overline{c^{\ve}}(t)\right)dt\nonumber\\
 &{+}\eta(0)\int_{\R^n}\nabla_x\varphi\cdot\ve\frac{1}{a+1}\,d\overline{vc_0^{\ve}}+\ve\frac{1}{a+1}\int_0^{\infty}\frac{d\eta}{dt}\int_{\R^n}\nabla_x\varphi\cdot\,d\overline{vc^{\ve}}(t)dt.\label{sc2w}
\end{align}
Further, we plug $\frac{d\eta}{dt}$ instead of $\eta$ for some $\eta\in C^2[0,\infty)$ into \cref{mom01w} and obtain that
\begin{align}
 &{-}\frac{d\eta}{dt}(0)\int_{\R^n}\varphi\,d\overline{c^{\ve}_0}-\int_0^{\infty}\frac{d^2\eta}{dt^2}\int_{\R^n}\varphi\,d\overline{c^{\ve}}(t)+\ve^{1-\kappa}\frac{d\eta}{dt}\int_{\R^n}\nabla_x\varphi\cdot\,d\overline{vc^{\ve}}(t)\,dt=0\label{mom0wt}
\end{align}
Combining \cref{mom01w}, \cref{sc2w}, and \cref{mom0wt} we arrive at \cref{mom012w}.
\end{proof}
\noindent
Now we can state and prove the following  upscaling {result}:
\begin{Theorem}[Parabolic limit]\label{convPL}
 Let  {\cref{AssqQ}} be satisfied. Assume that $\kappa=2$ and $\E[q]\equiv0$. For some $\ve_m\underset{m\rightarrow\infty}{\rightarrow}0$ let  $c^{\ve_m}_0\in {\cal M}_+(\R^n\times \overline{B}_{1}(0))$ be a sequence of initial data such that 
 \begin{align}
  c^{\ve_m}_0\underset{m\rightarrow\infty}{\overset{*}{\rightharpoonup}}c^{0}_{0}  \qquad\text{in }{\cal M}_+(\R^n\times \overline{B}_{1}(0))
 \end{align}
 for some $c^{0}_{0}\in {\cal M}_+(\R^n\times \overline{B}_{1}(0))$. Finally, let $c^{\ve_m}$ be a weak solution to \cref{transpce} in terms of {\cref{WKKTE}} corresponding to $c^{\ve_m}_0$. Then there exists a subsequence $c^{\ve_{m_k}}$ such that
 \begin{align}
  c^{\ve_{m_k}}\underset{m\rightarrow\infty}{\overset{*}{\rightharpoonup}}c^{0} \qquad\text{in }L^{\infty}_{w\text{-}*}(\R^+;{\cal M}_+(\R^n\times \overline{B}_{1}(0))),
 \end{align}
 where $c^{0}$ satisfies \cref{solc0} and $\overline{c^{0}}$ is a weak solution to \cref{PL} in terms of {\cref{WKPL}} corresponding to $\overline{c_0^0}$.
\end{Theorem}
\begin{Remark}
 Other than in \cite{HHW10} where the parabolic limit involved just self-diffusion of the cell population with a similar diffusion coefficient $\mathbb D[q]$, we do not require $q$ to be constant with respect to $x$. Such relaxation of conditions imposed on $q$ is highly relevant from the application viewpoint, since the orientational distribution  of tissue fibers usually varies from one point in space to the other. The spacial heterogeneity of the tissue structure plays a major role in haptotactic behaviour.
\end{Remark}

\begin{proof} (of \cref{convPL}) 
{To begin with, we observe that} since the sequence of initial measures is weak-$*$ converging, it is also uniformly bounded. Consequently, {the mass preservation property} \cref{conv} implies that  $c^{\ve_m}$ is uniformly bounded in $L^{\infty}_{w\text{-}*}(\R^+;{\cal M}_+(\R^n\times \overline{B}_{1}(0)))$. Due to the Banach-Alaoglu theorem (compare also {\cref{SecFSp}}), there exists a subsequence $c^{\ve_{m_k}}$ which weak-$*$ converges to some $c^{0}\in L^{\infty}_{w\text{-}*}(\R^+;{\cal M}_+(\R^n\times \overline{B}_{1}(0)))$. This allows to pass to the limit for $\ve=\ve_{m_k}$ as $k\rightarrow\infty$ in {equation \cref{mesoweaksol} from} \cref{WKKTE}. Thus, we obtain that
\begin{align}
 &\int_0^{\infty}\eta\int_{\R^n\times\overline{B}_{1}(0)}av\cdot\eta\varphi\nabla_v\psi\,dc^{0}(t)\,dt\nonumber\\
 =&\int_0^{\infty}\eta\left(n\int_{\R^n}\varphi\int_{\overline{B}_{1}(0)}\psi\,dq\,d\overline{c^{0}}(t)-\int_{\R^n\times\overline{B}_{1}(0)}\varphi\psi\,dc^{0}(t)\right)dt.
\end{align}
Resolving this distributional equation we conclude that
\begin{align}
 \int_{\overline{B}_{1}(0)}av\cdot\nabla_v\psi\,dc^{0}
 =&n\int_{\overline{B}_{1}(0)}\psi\,dq\,\overline{c^{0}}-\int_{\overline{B}_{1}(0)}\psi\,dc^{0}\qquad\text{in }L^{\infty}_{w\text{-}*}(\R^+;{\cal M}_+(\R^n)).\label{weakc0}
\end{align}
In particular, taking $\psi(v)=vv^T$ in \cref{weakc0} we {obtain} the expression \cref{M2} for  the second moment of $c^{0}$.
{Due to \cref{PropWMOME},  for each $m\in\N$ the weak solution $c^{\ve_m}$  satisfies the moment equation \cref{mom012w}.} Passing to the limit in  {this equation}, we obtain that
\begin{align}
 &{-}\eta(0)\int_{\R^n}\varphi\,d\overline{c^{0}_{0}}-\int_0^{\infty}\frac{d\eta}{dt}\int_{\R^n}\varphi\,d\overline{c^{0}}(t)\nonumber\\
 =&\int_0^{\infty}\eta\left(\int_{\R^n}\nabla_x\nabla_x^T\varphi: \frac{1}{a+1}\,d \overline{vv^Tc^{0}}(t)+\int_{\R^n}\nabla_x\varphi\cdot\frac{a}{a+1}\F\nabla_x Q\,d\overline{c^{0}}(t)\right)dt.\label{limpar}
\end{align}
Substituting \cref{M2} into \cref{limpar}, we arrive at the weak formulation {\cref{PLsol}} from \cref{WKPL}.

\noindent
For a fixed $\overline{c^{0}}$ a solution to \cref{weakc0} is given by \cref{solc0}. Since the equation is a linear one, this solution is also unique if the corresponding homogeneous equation
\begin{align}
 \int_{\overline{B}_{1}(0)}(av\cdot\nabla_v\psi+\psi)\,du
 =&0\qquad\text{for all }\psi\in C^1(\overline{B}_{1}(0))\label{hom}
\end{align}
has only the trivial solution in $\overline{B}_{1}(0)$. Observe that each  $g\in C_c(\overline{B}_{1}(0)\backslash\{0\})$ can be described as 
\begin{align}
 av\cdot\nabla_v\psi+\psi=g,
\end{align}
where 
\begin{align}
 \psi_g(v)=-|v|^{-\frac{1}{a}}\int_{|v|^{\frac{1}{a}}}^1{g}(y\sign(v))\,dy
\end{align}
belongs to $C^1(\overline{B}_{1}(0))$. Consequently, $\supp u\subset\{0\}$, i.e. $u=C\delta_0$ {where $\delta_0$ denotes the Dirac delta and  $C$ is some constant}. But then 
\begin{align}
 0=\int_{\overline{B}_{1}(0)}(av\cdot\nabla_v\psi+\psi)\,du=C\psi(0)\qquad\text{for all }\psi\in C^1(\overline{B}_{1}(0)),
\end{align}
which implies that $C=0$. This shows that $u\equiv0$. 
\end{proof}

An analogous result holds for the hyperbolic case:
\begin{Theorem}[Hyperbolic limit]\label{convHL}
 Let  {\cref{AssqQ}} be satisfied. Assume that $\kappa=1$. For some $\ve_m\underset{m\rightarrow\infty}{\rightarrow}0$ let  $c^{\ve_m}_0\in {\cal M}_+(\R^n\times \overline{B}_{1}(0))$ be a sequence of initial data such that 
 \begin{align}
  c^{\ve_m}_0\underset{m\rightarrow\infty}{\overset{*}{\rightharpoonup}}c^{0}_{0}  \qquad\text{in }{\cal M}_+(\R^n\times \overline{B}_{1}(0))
 \end{align}
 for some $c^{0}_{0}\in {\cal M}_+(\R^n\times \overline{B}_{1}(0))$. Finally, let $c^{\ve_m}$ be a weak solution to \cref{transpce} in terms of {\cref{WKKTE}} corresponding to $c^{\ve_m}_0$. Then there exists a subsequence $c^{\ve_{m_k}}$ such that
 \begin{align}
  c^{\ve_{m_k}}\underset{m\rightarrow\infty}{\overset{*}{\rightharpoonup}}c^{0} \qquad\text{in }L^{\infty}_{w\text{-}*}(\R^+;{\cal M}_+(\R^n\times \overline{B}_{1}(0))),\nonumber
 \end{align}
 where $c^{0}$ satisfies \cref{solc0} and $\overline{c^{0}}$ is a weak solution to \cref{HL} in terms of {\cref{WKHL}} corresponding to $\overline{c_0^0}$.
\end{Theorem}
We omit the proof of this Theorem since it is very similar to that of \cref{convPL}. 
\begin{Remark}
  In the same way as we have proved the rigorous convergence for the zero order approximations one could validate the other formal derivations performed in \cref{SecScaling}.
\end{Remark}

\section{Discussion and outlook}\label{SecDis}
 In recent years modelling with KTEs  in the 
 {multiscale modelling} framework has proved to be an effective approach to describing cell movement in a fibrous environment, as it carefully connects single cell dynamics with the evolution of {one or several} cell distribution function{s} depending on time, position, velocity, and possibly further activity variables.
Since such mesoscopic models are generally difficult to handle numerically, suitable  macroscopic approximations, such as, e.g. limits of parabolic or hyperbolic scalings, are often derived and solved instead. In this work we developed a new approach to dealing with such scalings for a general class of KTEs involving transport with respect to velocity. It relies on the method of characteristics and a  differential equation \cref{mom012} which connects moments of zero and second order. The latter key equation can be utilised for both parabolic and hyperbolic scalings  and offers a unified and transparent way of deriving macroscopic equations for approximations of an arbitrary high order. As an illustration, we  have deduced {DTE}s for zero and first order approximations on the macroscale for both mentioned scaling types for our KTE. It  turns out that our formal computations  can be  mimicked by the corresponding  operations with Radon measures.  Under rather general conditions on the parameters which allow for  a spatially heterogeneous measure-valued fiber orientation distribution, we have thus been able to validate our limit passages rigorously.

Most of the previous constructions  leading from KTEs to RDTEs rely on a  Hilbert {space} structure already on the level of formal scalings. Indeed, one typically assumes  the zero and first order approximations  to be orthogonal in a particular weighted space of square integrable functions. We have actually seen that this property fails to hold  for the model class considered here. 
Conversely, our approach does not rely on orthogonality in any way and is applicable to a broader class of KTEs.  

In addition to the macroscopic approximations of zero order or higher we have also  developed a transport equation which preserves positivity and the total mass and can be solved numerically in order to obtain a suitable mesoscopic first order  approximation  of the solutions  to the original KTE. 

The class of KTEs we have used here to illustrate our approach can account for a number of  motility features. On the macroscale,  it has led to {DTE}s which, depending on the chosen scaling type and approximation order,  include  such terms as: myopic diffusion, drift, and taxis with respect to a mesoscopic and/or macroscopic quantity. 
The latter macroscopic quantity  could be a hapto- or a chemoattractant. 
 In our model, the taxis with respect to such an attractant  is  caused by biochemical and/or biophysical stress perceived by the cells. It is modelled via Newton's second law in \cref{eq:microdyn} and includes flux-limitation. We have seen that flux-limited taxis can be recovered on the macroscale provided that an approximation of a sufficiently high order is used. 
 One way to extend our model would be to consider a more general form of acceleration in \cref{eq:microdyn}, e.g. by  letting the acceleration scaling coefficient depend on the attractant. A similar consideration was made in  \cite{Chauviere2007} regarding the so-called  'chemotaxis force'. 
  One could also consider dependencies on other macroscopic quantities,  including cell population density.  The latter, however, would require dealing with convergences in nonlinear terms and would render {a} rigorous limit passage considerably more difficult. This is because  the weak type of convergence used in this work  would no longer be sufficient in order to handle nonlinearities.
 Further, we have assumed the attractant to be some given function, thus allowing us to deal with a single equation accounting for the  cell motion. A more realistic model would have to include an equation characterising the dynamics of the attractant.  That  would be an ODE if  the tactic cue represents, for instance, volume fraction of tissue which is supposed to be degraded by tumour cells. If, on the other hand,  it represents the concentration of some chemoattractant, then we would need to consider a reaction-diffusion PDE with source terms characterising production by tumour cells depending on their local macroscopic density and  decay caused by other influences, along with a  linear diffusion. In both cases, however, we would then be dealing with a nonlinear, strongly coupled system. Once again this would make the rigorous analysis  much more challenging, if at all possible. 
  

All motility terms obtained in the macroscopic formulations carry some information about the underlying tissue structure: the drift and diffusion coefficients depend on the orientation distribution of tissue fibers, whereas  the chemotactic sensitivity tensor which controls the taxis with respect to a macroscopic attractant accounts for spacial heterogeneity. For simplicity we have taken the fiber distribution to be time-independent, assuming it to be some given function of spacial position and orientation. Relaxing this assumption would lead (for an example of a  formal, orthogonality-based deduction and numerical simulations of the obtained equations see \cite{CEKNSSW}) to a nonlinear strongly coupled meso-macro system with dynamically changing diffusion coefficient  and drift velocity which is highly challenging from the analytical point of view.

\phantomsection
\addcontentsline{toc}{section}{References}
\printbibliography


\appendix
\section{Appendix}\label{Append}
As previously announced in \cref{RemSolv}, in this final auxiliary Section we briefly touch on the solvability of the KTE \cref{meso} which we have upscaled in this work. 
\begin{Proposition}[Existence for the KTE]
Let \cref{AssqQ} and the assumptions of  \cref{WKKTE} (for $\ve=1$) be satisfied.
 Then there exists a  weak solution  to the KTE \cref{meso}.
\end{Proposition}
\begin{proof}{\it (Sketch)}
 To shorten the notation we introduce
\begin{align*}
 &z:=(x,v),\qquad
 V(t,z):=(v,S(t,z))^T.
\end{align*}
 Let us first assume that 
 in addition to conditions 1.-2. from \cref{AssqQ} it holds that
\begin{enumerate}
\setcounter{enumi}{2}
\item
$\nabla_x\nabla_x^T Q\in C_b(\R^+_0\times\R^n;\R^{n\times n})$, $\nabla_xF\in C_b(\R^n;\R^{n\times n\times n})$.
\end{enumerate}
Under the assumptions we made, standard ODE theory implies  that the ODE system  \cref{MicroD} is globally  uniquely solvable on $\R^n\times\R^n$. We denote  by $Z(t,s):\R^n\times\R^n\rightarrow \R^n\times\R^n$ the corresponding solution operator, meaning that $Z(t,s)(z_0)$ is the value  at time $t$ of the  solution of \cref{MicroD} which has started in $z_0$ at the initial time $s$. Again, the ODE theory implies that this map is well-defined and bijective for all $s,t\in\R$.
Moreover, exploiting the form of $S$, we see that
\begin{align}
 (Z(t,s))^{(-1)}(\R^n\times(\R^n\backslash {B}_{1}(0)))\subset \R^n\times(\R^n\backslash \overline{B}_{1}(0))\qquad\text{for }t>s.\label{nonexp}
\end{align}

Let us set
\begin{align}
 c_0:=0,\qquad q:=0\qquad\text{in }\R^n\times(\R^n\backslash \overline{B}_{1}(0)),\label{extens}
\end{align}
so that
\begin{align*}
 c_0\in {\cal M}_+(\R^n\times\R^n),\qquad q\in C_b(\R^n;{\cal M}_+(\R^n)).
\end{align*}
We recall that  the Cauchy problem for the conservative transport equation 
\begin{align}
 \partial_t \mu+\nabla_z\cdot(V\mu)=0\label{CL}
\end{align}
has a unique weak solution $\mu\in C_{w\text{-}*}(\R^+;{\cal M}_+(\R^n\times\R^n))$ for every initial value $\mu_0\in {\cal M}_+(\R^n\times\R^n)$, see e.g. \cite[Chapter 5 Theorem 5.34]{Villani}. Moreover, its solution is given by an explicit formula:
\begin{align}
 \mu(t)=Z(t,0)\#\mu_0,\label{solCL}
\end{align}
where $T\#\mu$ denotes the push-forward of measure $\mu$ under a map $T:\R^n\times\R^n\rightarrow\R^n\times\R^n$.  We recall that for $T$  bijective the total mass is preserved under this operator:
\begin{align}
 \int_{\R^n\times\R^n}\,d(T\#\mu_0)=\int_{\R^n\times\R^n}\,d\mu_0.\label{masspres}
\end{align}
 This is the  case then for $T=Z(t,s)$ for any $s,t\in\R$.

Turning to the KTE \cref{meso}, we rewrite this equation in the following form: for all $t\in\R^+_0$
\begin{align}
 c(t)=Z(t,0)\#c_0+\int_0^tZ(t,s)\#(nq\overline{c}-c)(s)\,ds\qquad\text{in }{\cal M}_+(\R^n\times\R^n).\label{FPP}
\end{align}
Here and below the integration is understood in the weak-$*$ sense. 
Using the Banach fixed-point theorem and then a standard extension argument, one readily verifies that equation \cref{FPP} is solvable in $C_{w\text{-}*}(\R^+;{\cal M}_+(\R^n\times\R^n))$, and that its solution is unique on every time interval. Using \cref{masspres} and the assumptions on $q$, it is straightforward to check that 
 the total mass  is preserved over time:
 \begin{align}
 \int_{\R^n\times\R^n}\,dc(t)=\int_{\R^n\times\R^n}\,dc_0.\label{mpres}
\end{align}


Next, we verify that solutions to \cref{FPP} are weak solutions to the KTE \cref{meso}. 
Using \cref{CL,solCL,FPP}, we compute: for all $t\in\R^+_0$
\begin{align}
 \partial_tc(t)=&\partial_t (Z(t,0)\#c_0)+\int_0^t\partial_t(Z(t,s)\#(nq\overline{c}-c)(s))\,ds+(nq\overline{c}-c)(t)\nonumber\\
 =&-\nabla_z\cdot (V(t,\cdot)Z(t,0)\#c_0)-\int_0^t\nabla_z\cdot(V(t,\cdot)Z(t,s)\#(nq\overline{c}-c)(s))\,ds+(nq\overline{c}-c)(t)\nonumber\\
 =&-\nabla_z\cdot\left(V(t,\cdot)\left(Z(t,0)\#c_0+\int_0^tZ(t,s)\#(nq\overline{c}-c)(s)\,ds\right)\right)+(nq\overline{c}-c)(t)\nonumber\\
 =&-\nabla_z\cdot\left(V(t,\cdot)c(t)\right)+(nq\overline{c}-c)(t)\qquad\text{in }(C^1_c(\R^n\times\R^n))^*,\label{weakeq}
\end{align}
as  required. 
Finally, thanks to \cref{extens,nonexp} a weak solution  to \cref{weakeq} satisfies for all $t>0$
\begin{align}
 c(t)=e^{-t}Z(t,0)\#c_0=0\qquad \text{in }\R^n\times(\R^n\backslash {B}_{1}(0)).\label{vanish}
\end{align}
Combining \cref{weakeq,vanish,mpres}, we conclude that the solution to \cref{FPP} is a weak solution to the KTE \cref{meso} in terms of \cref{WKKTE}. A standard approximation argument finally yields existence under the original \cref{AssqQ}.
\end{proof}
\begin{Remark}[Uniqueness]
 Under some additional smoothness of $Q$ and $q$ with respect to $x$ variable one can verify the uniqueness of weak solutions to \cref{meso}.
\end{Remark}

\end{document}